*Research Article*

# The Implicit Midpoint Procedures for Asymptotically Nonexpansive Mappings

M. O. Aibinu ,[1,2] S. C. Thakur,[2] and S. Moyo[3]

[1]*Institute for Systems Science, Durban University of Technology, Durban 4000, South Africa*
[2]*KZN CoLab, Durban University of Technology, Durban 4000, South Africa*
[3]*Institute for Systems Science & Office of the DVC Research, Innovation & Engagement Milena Court,
 Durban University of Technology, Durban 4000, South Africa*

Correspondence should be addressed to M. O. Aibinu; moaibinu@yahoo.com

The concept of asymptotically nonexpansive mappings is an important generalization of the class of nonexpansive mappings. Implicit midpoint procedures are extremely fundamental for solving equations involving nonlinear operators. This paper studies the convergence analysis of the class of asymptotically nonexpansive mappings by the implicit midpoint iterative procedures. The necessary conditions for the convergence of the class of asymptotically nonexpansive mappings are established, by using a well-known iterative algorithm which plays important roles in the computation of fixed points of nonlinear mappings. A numerical example is presented to illustrate the convergence result. Under relaxed conditions on the parameters, some algorithms and strong convergence results were derived to obtain some results in the literature as corollaries.

## 1. Introduction

Let $K$ be a bounded subset of a Banach space $E$. A mapping $T: K \longrightarrow K$ is called asymptotically nonexpansive if there exists a sequence $\{k_n\}$ of positive real numbers with $k_n \longrightarrow 1$ as $n \longrightarrow \infty$ for which

$$\|T^n x - T^n y\| \leq k_n \|x - y\|, \quad \text{for all } x, y \in K \text{ and } n \in \mathbb{N}. \tag{1}$$

In 1972, Geobel and Kirk [1] introduced this concept of asymptotically nonexpansive mappings as an important generalization of the class of nonexpansive mappings in which $k_n = 1$ for all $n$ [1]. $T$ is said to be uniformly $L$-Lipschitzian, if there exists a constant $L > 0$ such that

$$\|T^n x - T^n y\| \leq L\|x - y\|, \quad \text{for all } x, y \in K \text{ and } n \in \mathbb{N}. \tag{2}$$

Notice that every asymptotically nonexpansive is uniformly $L$-Lipschitzian with a constant $L = \sup_{n \in \mathbb{N}} k_n \geq 1$. The set of fixed point of $T, \{x \in K : Tx = x\}$, will be denoted by $F(T)$. A mapping $f: K \longrightarrow K$ is called a $\alpha$-contraction if there exists $\alpha \in [0, 1)$ such that

$$\|f(x) - f(y)\| \leq \alpha \|x - y\|, \quad \text{for all } x, y \in K. \tag{3}$$

It is well known that a contraction $f$ on $K$ has a unique fixed point in $K$.

A powerful numerical method for solving ordinary differential equations and differential algebraic equations, which also has a long history, is the implicit midpoint procedure. Akin to the implicit midpoint procedure is the fractal structures of Mandelbrot and Julia sets which have been demonstrated for practical application in quadratic, cubic, and higher degree polynomials. Kang et al. in [2] defined Jungck Noor iteration with $s$-convexity and established the escape criterions for quadratic, cubic, and $n$th degree complex polynomials. Auzinger and Frank in [3] studied the structure of the global discretization error for the implicit midpoint and trapezoidal rules applied to nonlinear stiff initial value problems. Bader and Deuflhard in [4] introduced a semi-implicit extrapolation method especially designed for the numerical solution of stiff systems of ordinary differential equations. The implicit midpoint rule is described as a theoretical foundation of the numerical treatment of problems arising in physical and biological sciences (Kastner-Maresch [5]). The implicit midpoint rule is applied to obtain the periodic solution of a nonlinear time-dependent evolution equation and a Fredholm integral equation (Somali and Davulcua [6]). Consider the ordinary differential equation:

$$\begin{aligned} x' &= f(t), \\ x(0) &= x_0. \end{aligned} \quad (4)$$

A sequence $\{x_n\}$ is generated by the implicit midpoint rule via the recursion:

$$\frac{1}{h}(x_{n+1} - x_n) = f\left(\frac{x_n + x_{n+1}}{2}\right), \quad n \in \mathbb{N}, \quad (5)$$

where $h > 0$ is a stepsize and $\mathbb{N}$ is the set of positive integers. For a Lipschitz continuous and sufficiently smooth map $f: \mathbb{R}^N \longrightarrow \mathbb{R}^N$, it is known that, as $h \longrightarrow 0$ uniformly over $t \in [0, \overline{t}]$ for any fixed $\overline{t} > 0$, the sequence $\{x_n\}$ generated by (5) converges to the exact solution of (4). By rewriting the function $f$ in the form $f(t) = g(t) - t$, the differential equation (4) becomes $x\prime = g(t) - t$. Consequently, the equilibrium problem associated with the differential equation is transformed to the fixed point problem $g(t) = t$ [7]. Extension of the implicit midpoint rule to nonexpansive mappings by Alghamdi et al. in [8] generates a recursion sequence:

$$x_{n+1} = (1 - a_n)x_n + a_n T\left(\frac{x_n + x_{n+1}}{2}\right), \quad n \in \mathbb{N}, \quad (6)$$

where $a_n \in (0, 1)$, $T: K \longrightarrow K$ is a nonexpansive mapping and $K$ is a closed convex subset of a real Hilbert space $H$. They proved the weak convergence of (6) under certain conditions on $\{a_n\}$. Still, in Hilbert space, Xu et al. in [9] used contractions to regularize the implicit midpoint rule (6) and introduced the implicit procedure:

$$x_{n+1} = a_n f(x_n) + (1 - a_n) T\left(\frac{x_n + x_{n+1}}{2}\right), \quad n \in \mathbb{N}. \quad (7)$$

They proved a strong convergence theorem for the sequence $\{x_n\}$ to a fixed point $p$ of $T$ which also solves the variational inequality:

$$\langle (I - f)p, x - p \rangle \geq 0, \quad \forall x \in F(T). \quad (8)$$

In 2015, Yao et al. in [7] introduced

$$x_{n+1} = a_n f(x_n) + b_n x_n + c_n T\left(\frac{x_n + x_{n+1}}{2}\right), \quad n \in \mathbb{N}, \quad (9)$$

which gives a faster approximation compared with (7), where $T$ is a nonexpansive mapping in a Hilbert space and $a_n + b_n + c_n = 1, \forall n \in \mathbb{N}$. Aibinu et al. in [10] compared the rate of convergence of the iteration procedures (6), (7), and (9). Aibinu and Kim in [11] recently introduced a new scheme of the viscosity implicit iterative algorithms for nonexpansive mappings in Banach spaces. Suitable conditions were imposed on the control parameters to prove a strong convergence theorem for the considered iterative sequence. Also, Aibinu and Kim in [12] studied the analytical comparison of schemes (7) and (9) to determine the sequence that converges faster in approximating a fixed point of a nonexpansive mapping. Literature review reveals the following problems, which this paper is devoted to address.

*Problem 1.* Can one study the implicit iterative procedure (9) for the class of asymptotically nonexpansive mappings which is more general than nonexpansive mappings?

*Problem 2.* Under what conditions will the main results of Yao et al. in [7] hold for asymptotically nonexpansive mappings in the general Banach spaces?

Motivated by the previous works, this paper is devoted for the extension of the previous results in the literature, to a more general space and to study the analogue of algorithm (9) for the class of asymptotically nonexpansive mappings. Therefore, for a Banach space with a uniformly Gâteaux differentiable norm possessing uniform normal structure and for arbitrary $x_1 \in K$, this paper considers the iterative algorithm given by

$$x_{n+1} = a_n f(x_n) + b_n x_n + c_n T^n\left(\frac{x_n + x_{n+1}}{2}\right), \quad n \in \mathbb{N}, \quad (10)$$

where $T$ is an asymptotically nonexpansive mapping and the real sequences $\{a_n\} \subset (0, 1)$, $\{b_n\} \subset [0, 1)$ and $\{c_n\} \subset (0, 1)$ are chosen such that $a_n + b_n + c_n = 1 \forall n \in \mathbb{N}$. We study the convergence of the important class of asymptotically nonexpansive mappings, which is a generalization of the class of nonexpansive mappings. The necessary conditions for the convergence of the class of asymptotically nonexpansive mappings are established, by using this well-known iterative algorithm which plays important roles in the computation of fixed points of nonlinear mappings. A numerical example is given to illustrate the convergence result.

## 2. Preliminaries

Let $E$ be a real Banach space with dual $E^*$ and denote the norm on $E$ by $\|.\|$. The normalized duality mapping $J: E \longrightarrow 2^{E^*}$ is defined as

$$J(x) = \{f \in E^*: \langle x, f \rangle = \|x\| \|f\|, \|x\| = \|f\|\}, \quad (11)$$

where $\langle ., . \rangle$ is the duality pairing between $E$ and $E^*$. Let $B_E$ denotes the unit ball of $E$. The modulus of convexity of $E$ is defined as follows $\delta_E(\varepsilon) = \inf\{1 - (\|x + y\|/2): x, y \in B_E, \|x - y\| \geq \varepsilon\}, 0 \leq \varepsilon \leq 2$. $E$ is uniformly convex if and only if $\delta_E(\varepsilon) > 0$ for every $\varepsilon \in (0, 2]$. $E$ is said to be smooth (or Gâteaux differentiable) if the limit,

$$\lim_{t \longrightarrow 0^+} \frac{\|x + ty\| - \|x\|}{t}, \quad (12)$$

exists for each $x, y \in B_E$. $E$ is said to have uniformly Gâteaux differentiable norm if, for each $y \in B_E$, the limit is attained uniformly for $x \in B_E$ and uniformly smooth if it is smooth and the limit is attained uniformly for each $x, y \in B_E$. Recall that if $E$ is smooth, then $J$ is single-valued and onto if $E$ is reflexive. Furthermore, the normalized duality mapping $J$ is uniformly continuous on bounded subsets of $E$ from the strong topology of $E$ to the weak-star topology of $E^*$ if $E$ is a Banach space with a uniformly Gâteaux differentiable norm.

Let $K$ be a nonempty bounded closed convex subset of a Banach space $E$, and let the diameter of $K$ be defined by $d(K) = \sup\{\|x - y\|: x, y \in K\}$. For each $x \in K$, let $r(x, K) = \sup\{\|x - y\|: y \in K\}$ and let $r(K) = \inf(r(x, K): x \in K)$, and the Chebyshev radius of $K$ is relative to itself. The normal structure coefficient of $E$ (see [13]) is defined as the number:

$$N(E) = \inf\left\{\frac{d(K)}{r(K)}: K \text{ bounded closed convex subset of } E \text{ with } d(K) > 0\right\}. \tag{13}$$

A space $E$, such that $N(E) > 1$, is said to have uniformly normal structure. Recall that a space with uniformly normal structure is reflexive and that all uniformly convex or uniformly smooth Banach spaces have uniform normal structure [14, 15]. If $E$ is a reflexive Banach space with modulus of convexity $\delta_E$, then $N(E) \geq (1 - \delta_E(1))^{-1}$. It was proved that if the space $E$ is uniformly convex; then, every asymptotically nonexpansive self-mapping $T$ of $K$ has a fixed point [1]. Also, it has been proved that if $\liminf \rho(t)/t < 1/2$, then $N(E) > 1$, where $\rho(t) = \sup\{1/2(\|x + ty\| + \|x - ty\|) - 1, x, y \in B_E\}$ for $t \geq 0$ [16, 17].

Recall that a function $f: X \longrightarrow R \cup \{+\infty\}$ is said to be weakly lower semicontinuous at $x_0 \in X$ if whenever $\{x_n\}$ is a sequence in $X$ such that $x_n$ converges weakly to $x_0$; then,

$$f(x_0) \leq \liminf_{n \longrightarrow \infty} f(x_n). \tag{14}$$

According to Jung and Kim in [18], let $\mu$ be a mean on positive integers $\mathbb{N}$. A mean $\mu$ on $\mathbb{N}$ is a continuous linear functional on $l^\infty$ satisfying $\|\mu\| = 1 = \mu(1)$. It is called a Banach limit if $\mu_n(a_n) = \mu_n(a_{n+1})$ for every $a = (a_1, a_2, \ldots) \in l^\infty$.

The following lemmas are needed in the sequel.

**Lemma 1** (see [19]). *Let $K$ be a nonempty closed convex subset of a Banach space $E$ with a uniformly Gâteaux differentiable norm, let $\{x_n\}$ be a bounded sequence of $E$, and let $\mu_n$ be a mean on $\mathbb{N}$. Let $x \in K$. Then,*

$$\mu_n \|x_n - x\|^2 = \min_{z \in K} \mu_n \|x_n - z\|^2, \tag{15}$$

*if and only if $\mu_n \langle z - x, J(x_n - x) \rangle \leq 0$ for all $z \in K$, where $J$ is the duality mapping of $E$.*

**Lemma 2** (see [20]). *Let $\{\sigma_n\}$ be a sequence of nonnegative real numbers satisfying the property:*

$$\sigma_{n+1} = (1 - \gamma_n)\sigma_n + \gamma_n \beta_n, \quad n > 0, \tag{16}$$

*where $\{\gamma_n\} \subset (0, 1)$ and $\{\beta_n\} \subset \mathbb{R}$ such that*

(i) $\sum_{n=1}^\infty \gamma_n = \infty$

(ii) $\limsup_{n \longrightarrow \infty} \beta_n \leq 0$

*Then, $\{\sigma_n\}$ converges to zero, as $n \longrightarrow \infty$.*

**Lemma 3** (see [21]). *Let $E$ be a Banach space with uniform normal structure, $K$ a nonempty closed convex and bounded subset of $E$, and $T: K \longrightarrow K$ an asymptotically nonexpansive mapping. Then, $T$ has a fixed point.*

**Lemma 4** (see [22]). *Let $E$ be a real Banach space and $J$ the normalized duality map on $E$. Then, for any given $x, y \in E$, the following inequality holds:*

$$\|x + y\|^2 \leq \|x\|^2 + 2\langle y, j(x + y) \rangle, \quad \forall j(x + y) \in J(x + y). \tag{17}$$

**Lemma 5** (see [23, 24]). *Let $\eta$ be a real number and $(x_0, x_1, \ldots) \in l^\infty$ such that $\mu_n x_n \leq \eta$ for all Banach limits. If $\limsup_{n \longrightarrow \infty}(x_{n+1} - x_n) \leq 0$, then $\limsup_{n \longrightarrow \infty} x_n \leq \eta$.*

**Lemma 6** (see [15]). *Suppose $E$ is a Banach space with uniformly normal structure, $K$ is a nonempty bounded subset of $E$, and $T: K \longrightarrow K$ is a uniformly $k$-Lipschitzian mapping with $k < N(E)^{1/2}$. Suppose also there exists a nonempty, bounded, closed, and convex subset $G$ of $K$ with the following property:*

$$(P): p \in G \text{ implies } \omega_w(p) \in G, \tag{18}$$

*where $\omega_w(p)$ is the weak $\omega$-limit set of $T$ at $p$, that is, the set*

$$\left\{x \in E: x = \text{weak} - \lim_{j \longrightarrow \infty} T^{n_j} p \text{ for some } n_j \longrightarrow \infty\right\}. \tag{19}$$

*Then, $T$ has a fixed point in $G$.*

## 3. Main Results

*Assumption 1.* Let $E$ be a real Banach space with a uniformly Gâteaux differentiable norm possessing uniform normal structure and $K$ a nonempty bounded closed convex subset of $E$. Suppose $T: K \longrightarrow K$ is an asymptotically nonexpansive mapping with $F(T) \neq \emptyset$ and $f: K \longrightarrow K$ is a contraction with constant $\alpha \in [0, 1)$. For arbitrary $x_1 \in K$ and real sequences $\{a_n\} \subset (0, 1), \{b_n\} \subset [0, 1)$ and $\{c_n\} \subset (0, 1)$ satisfying $a_n + b_n + c_n = 1, \forall n \in \mathbb{N}$, this paper considers the iterative scheme $\{x_n\}$ given by (10). Assume that $\{a_n\}$ satisfies the following conditions:

(i) $\lim_{n \to \infty} a_n = 0$

(ii) $\sum_{n=1}^{\infty} a_n = \infty$

(iii) $\lim_{n \to \infty} ((k_n^2 - 1)/a_n) = 0$

where $\{k_n\}$ is a sequence of positive real numbers with $k_n \to 1$ as $n \to \infty$. Then, the sequence $\{x_n\}$ is well defined. To show this, for arbitrary $\omega \in K$, define the mapping $T_\omega: K \to K$:

$$x \mapsto T_\omega := a_n f(\omega) + b_n \omega + c_n T^n \left( \frac{\omega + x}{2} \right), \quad n \in \mathbb{N}. \quad (20)$$

Then, for $x, y \in K$,

$$\|T_\omega x - T_\omega y\| \leq c_n \left( T^n \left( \frac{\omega + x}{2} \right) - T^n \left( \frac{\omega + y}{2} \right) \right) \leq \frac{c_n k_n}{2} \|x - y\|. \quad (21)$$

Using condition (iii), for any given positive number $\varepsilon, 0 < \varepsilon < 1 - \alpha$, there exists a sufficient large positive integer $n_0$, such that, for any $n \geq n_0$, we obtain that

$$k_n^2 - 1 \leq 2\varepsilon a_n,$$

$$k_n - 1 \leq \frac{k_n + 1}{2}(k_n - 1) \leq \frac{k_n^2 - 1}{2} \leq \varepsilon a_n. \quad (22)$$

Therefore, since $\{k_n\}$ is a sequence of positive real numbers with $k_n \to 1$ as $n \to \infty$ and $k_n - 1 \leq \varepsilon a_n$ for all $n \geq n_0$, it is obvious that

$$\begin{aligned}
c_n k_n &= 1 + (k_n - 1) - k_n + c_n k_n \\
&\leq 1 + \varepsilon a_n - k_n(1 - c_n) \\
&= 1 + \varepsilon a_n - (1 - c_n) \\
&= 1 + \varepsilon a_n - (a_n + b_n) \\
&= 1 - a_n(1 - \varepsilon) - b_n \\
&\leq 1.
\end{aligned} \quad (23)$$

Thus, from (21) and (23), it shows that $T_\omega$ is a contraction. Therefore, (10) is well defined since every contraction in a Banach space has a fixed point.

Firstly, the proof of the following lemmas are given, which are useful in establishing the main result.

**Lemma 7.** *Let $E$ be a Banach space with a uniform normal structure, $K$ a nonempty bounded closed convex subset of $E$, and $T: K \to K$ an asymptotically nonexpansive mapping. Suppose $f: K \to K$ is a contraction with constant $\alpha \in [0, 1)$ and assume that $\{a_n\}$ satisfies conditions (i) – (iii). For an arbitrary $x_1 \in K$, the iterative sequence $\{x_n\}$ given by (10) is bounded.*

*Proof.* The sequence $\{x_n\}$ is shown to be bounded.

By Lemma 3, $F(T) \neq \emptyset$. Then, for $p \in F(T)$,

$$\begin{aligned}
\|x_{n+1} - p\| &= \left\| a_n f(x_n) + b_n x_n + c_n T^n \left( \frac{x_n + x_{n+1}}{2} \right) - p \right\| \\
&= \left\| a_n (f(x_n) - f(p)) + a_n (f(p) - p) + b_n (x_n - p) \right. \\
&\quad \left. + c_n \left( T^n \left( \frac{x_n + x_{n+1}}{2} \right) - p \right) \right\| \\
&\leq a_n \|f(x_n) - f(p)\| + a_n \|f(p) - p\| + b_n \|x_n - p\| \\
&\quad + c_n \left\| T^n \left( \frac{x_n + x_{n+1}}{2} \right) - p \right\| \\
&\leq a_n \|f(x_n) - f(p)\| + a_n \|f(p) - p\| + b_n \|x_n - p\| \\
&\quad + c_n k_n \left\| \frac{x_n + x_{n+1}}{2} - p \right\| \\
&\leq \alpha a_n \|x_n - p\| + a_n \|f(p) - p\| + b_n \|x_n - p\| + \frac{c_n k_n}{2} \|x_n - p\| \\
&\quad + \frac{c_n k_n}{2} \|x_{n+1} - p\|.
\end{aligned} \quad (24)$$

Therefore,

$$\left(1 - \frac{c_n k_n}{2}\right) \|x_{n+1} - p\| \leq \left(\alpha a_n + b_n + \frac{c_n k_n}{2}\right) \|x_n - p\| + a_n \|f(p) - p\|, \quad (25)$$

which is equivalent to

$$\frac{2 - c_n k_n}{2} \|x_{n+1} - p\| \leq \frac{2\alpha a_n + 2b_n + c_n k_n}{2} \|x_n - p\| + a_n \|f(p) - p\|. \quad (26)$$

Observe that

$$\begin{aligned}
2 - c_n k_n &= 2 - (1 - (a_n + b_n))k_n \\
&= 2 - (k_n - k_n(a_n + b_n)) \\
&= 1 - (k_n - 1) + k_n(a_n + b_n) \\
&\geq 1 - \varepsilon a_n + k_n(a_n + b_n) \\
&= 1 + (k_n - \varepsilon)a_n + k_n b_n \\
&\geq 1 + (1 - \varepsilon)a_n + b_n,
\end{aligned} \quad (27)$$

for sufficiently large $n$.

Also,

$$\begin{aligned}
2\alpha a_n + 2b_n + c_n k_n &= 2\alpha a_n + 2b_n + (1 - (a_n + b_n))k_n \\
&= 2\alpha a_n + 2b_n + k_n - k_n(a_n + b_n) \\
&= 1 + (k_n - 1) + 2\alpha a_n - k_n a_n + 2b_n - k_n b_n \\
&\leq 1 + \varepsilon a_n + 2\alpha a_n - k_n a_n + 2b_n - k_n b_n \\
&= 1 - (k_n - 2\alpha - \varepsilon)a_n - (k_n - 2)b_n \\
&\leq 1 - (1 - 2\alpha - \varepsilon)a_n - (1 - 2)b_n \\
&= 1 - (1 - 2\alpha - \varepsilon)a_n + b_n,
\end{aligned}$$

for sufficiently large $n$.

$$(28)$$

Consequently, (26) becomes

$$\|x_{n+1} - p\| \leq \frac{1 - (1 - 2\alpha - \varepsilon)a_n + b_n}{1 + (1 - \varepsilon)a_n + b_n}\|x_n - p\| + \frac{2a_n}{1 + (1 - \varepsilon)a_n + b_n}\|f(p) - p\|$$

$$= \left(1 - \frac{2(1 - \alpha - \varepsilon)a_n}{1 + (1 - \varepsilon)a_n + b_n}\right)\|x_n - p\| + \frac{2(1 - \alpha - \varepsilon)a_n}{1 + (1 - \varepsilon)a_n + b_n}\left(\frac{1}{1 - \alpha - \varepsilon}\|f(p) - p\|\right)$$

$$\leq \max\left\{\|x_n - p\|, \frac{1}{1 - \alpha - \varepsilon}\|f(p) - p\|\right\} \quad (29)$$

$$\vdots$$

$$\leq \max\left\{\|x_1 - p\|, \frac{1}{1 - \alpha - \varepsilon}\|f(p) - p\|\right\}.$$

This implies that the sequence $\{x_n\}$ is bounded and hence $\{f(x_n)\}$ and $\{T^n x_n\}$ are also bounded.

Obviously,

$$\|T^n x_n\| = \|T^n x_n - p + p\|$$
$$\leq \|T^n x_n - p\| + \|p\|$$
$$\leq k_n \|x_n - p\| + \|p\|$$
$$\vdots$$
$$\leq R \max\left\{\|x_1 - p\|, \frac{1}{1 - \alpha - \varepsilon}\|f(p) - p\|\right\} + \|p\|, \quad (30)$$

where $R := \sup k_n$. □

**Lemma 8.** *Let $E$ be a real Banach space with a uniformly Gateaux differentiable norm possessing uniform normal structure, $K$ a nonempty bounded closed convex subset of $E$, $T: K \longrightarrow K$ an asymptotically nonexpansive mapping with sequence $\{k_n\} \subset [1, \infty)$, $f: K \longrightarrow K$ a continuous mapping, and $\{x_n\}$ a bounded sequence in $K$ such that $\lim_{n\to\infty}\|x_n - T^n x_n\| = 0$. Suppose $\{z_t\}$ is a path in $K$ defined by $z_t = tf(z_t) + (1-t)T^n z_t, t \in (0,1)$ such that $z_t \longrightarrow z$ as $t \longrightarrow 0^+$. Then,*

$$\limsup_{n\to\infty} \langle f(z) - z, J(x_n - z)\rangle \leq 0. \quad (31)$$

*Proof.* From $\lim_{n\to\infty}\|x_n - T^n x_n\| = 0$, let $M := \sup_{n\to\infty}\|x_n - T^n x_n\|$ and define $\{t_n\} := \{\sqrt{\|x_n - T^n x_n\|}/M\}, n \in \mathbb{N}$. Clearly,

$$\lim_{n\to\infty} \frac{\|x_n - T^n x_n\|}{t_n} = 0. \quad (32)$$

Since $z_t$ is given by $z_t = tf(z_t) + (1-t)T^n z_t$, we can write

$$z_{t_n} - x_n = t_n(f(z_{t_n}) - x_n) + (1 - t_n)(T^n z_{t_n} - x_n). \quad (33)$$

Since $\{x_n\}, \{z_{t_n}\}$, and $\{T^n x_n\}$ are bounded, $\{z_{t_n} - x_n\}$ is bounded. By Lemma 4,

$$\|z_{t_n} - x_n\|^2 \leq (1 - t_n)^2 \|T^n z_{t_n} - x_n\|^2 + 2t_n \langle f(z_{t_n}) - x_n, J(z_{t_n} - x_n)\rangle$$

$$\leq (1 - t_n)^2 \left(\|T^n z_{t_n} - T^n x_n\| + \|T^n x_n - x_n\|\right)^2$$

$$+ 2t_n \langle f(z_{t_n}) - z_{t_n}, J(z_{t_n} - x_n)\rangle + 2t_n \|z_{t_n} - x_n\|^2$$

$$\leq (1 - t_n)^2 \left(k_n \|z_{t_n} - x_n\| + \|x_n - T^n x_n\|\right)^2$$

$$+ 2t_n \langle f(z_{t_n}) - z_{t_n}, J(z_{t_n} - x_n)\rangle + 2t_n \|z_{t_n} - x_n\|^2$$

$$\leq (1 - t_n)^2 \left(k_n^2 \|z_{t_n} - x_n\|^2 + \|x_n - T^n x_n\| \left(2k_n \|z_{t_n} - x_n\| + \|x_n - T^n x_n\|\right)\right)$$

$$+ 2t_n \langle f(z_{t_n}) - z_{t_n}, J(z_{t_n} - x_n)\rangle + 2t_n \left\{\|z_{t_n} - x_n\|\right\}^2. \quad (34)$$

Therefore,

$$\langle f(z_{t_n}) - z_{t_n}, J(x_n - z_{t_n})\rangle \leq \frac{k_n^2(1-t_n)^2 - (1-2t_n)}{2t_n}\|z_{t_n} - x_n\|^2$$

$$+ \frac{(1-t_n)^2\|x_n - T^n x_n\|}{2t_n}\left(2k_n\|z_{t_n} - x_n\| + \|x_n - T^n x_n\|\right)$$

$$\leq \frac{(k_n^2 - 1)(1-2t_n) + k_n^2 t_n^2}{2t_n}\|z_{t_n} - x_n\|^2$$

$$+ \frac{(1+t_n^2)\|x_n - T^n x_n\|}{2t_n}\left(2k_n\|z_{t_n} - x_n\| + \|x_n - T^n x_n\|\right). \tag{35}$$

From (32) and by recalling that $k_n \longrightarrow 1$ as $n \longrightarrow \infty$,

$$\limsup_{n\longrightarrow\infty}\langle f(z_{t_n}) - z_{t_n}, J(x_n - z_{t_n})\rangle \leq 0. \tag{36}$$

Observe that

$$\langle f(z) - z, J(x_n - z)\rangle = \langle f(z) - z, J(x_n - z) - J(x_n - z_{t_n})\rangle$$
$$+ \langle f(z) - f(z_{t_n}), J(x_n - z_{t_n})\rangle + \langle z_{t_n} - z, J(x_n - z_{t_n})\rangle$$
$$+ \langle f(z_{t_n}) - z_{t_n}, J(x_n - z_{t_n})\rangle. \tag{37}$$

Noting the hypothesis that $z_{t_n} \longrightarrow z$ as $n \longrightarrow \infty$, coupled with the continuity of $f$ and uniform continuity of the duality mapping $J$ on bounded subsets of $E$, it gives

$$\langle z_{t_n} - z, J(x_n - z_{t_n})\rangle \longrightarrow 0, \quad n \longrightarrow \infty,$$
$$\langle f(z) - f(z_{t_n}), J(x_n - z_{t_n})\rangle \longrightarrow 0, \quad n \longrightarrow \infty,$$
$$\langle f(z) - z, J(x_n - z) - J(x_n - z_{t_n})\rangle \longrightarrow 0, \quad n \longrightarrow \infty. \tag{38}$$

From (36)–(38), it can be deduced that

$$\limsup_{n\longrightarrow\infty}\langle f(z) - z, J(x_n - z)\rangle \leq 0. \tag{39}$$
□

**Theorem 1.** *Let $E$ be a real Banach space with a uniformly Gateaux differentiable norm possessing uniform normal structure, $K$ a nonempty bounded closed convex subset of $E$, and $T: K \longrightarrow K$ an asymptotically nonexpansive mapping with the sequence $\{k_n\}$ of positive real numbers such that $k_n \longrightarrow 1$ as $n \longrightarrow \infty$, $\sup_{n\in\mathbb{N}} k_n \leq N(E)^{1/2}$, and $F(T) \neq \emptyset$. Suppose $f: K \longrightarrow K$ is a contraction with constant $\alpha \in [0, 1)$. Assume that $\{a_n\}$ satisfies conditions (i) – (iii) and in addition (iv) $\lim_{n\longrightarrow\infty}\|x_n - T^n x_n\| = 0$*

*Then, for an arbitrary $x_1 \in K$, the iterative sequence $\{x_n\}$ given by (10) converges in norm to a fixed point $p$ of $T$ which is the unique solution of variational inequality:*

$$\langle (I - f)p, j(x - p)\rangle \geq 0, \quad \forall x \in F(T). \tag{40}$$

*Proof.* □

*Step 1.* It is shown in this step that $\lim_{n\longrightarrow\infty}\|x_n - Tx_n\| = 0$.
Observe that

$$\|x_n - Tx_n\| = \|x_n - x_{n+1} + x_{n+1} - T^n x_n + T^n x_n - Tx_n\|$$
$$\leq \|x_n - x_{n+1}\| + \|x_{n+1} - T^n x_n\| + \|T^n x_n - Tx_n\|$$
$$\leq \|x_{n+1} - x_n\| + \|x_{n+1} - T^n x_n\| + k_1 \|x_n - T^{n-1} x_n\|, \tag{41}$$

for some $k_1 \in \mathbb{R}^+$, where $\mathbb{R}^+$ denotes the set of positive real numbers. It is required that

$$\lim_{n\longrightarrow\infty}\|x_{n+1} - x_n\| = \lim_{n\longrightarrow\infty}\|x_{n+1} - T^n x_n\| = \lim_{n\longrightarrow\infty}\|x_n - T^{n-1} x_n\| = 0. \tag{42}$$

To evaluate $\lim_{n\longrightarrow\infty}\|x_{n+1} - x_n\|$,

$$\|x_{n+1} - x_n\| \leq \|x_{n+1} - T^n x_n\| + \|T^n x_n - x_n\|$$

$$= \left\|a_n f(x_n) + b_n x_n + c_n T^n\left(\frac{x_n + x_{n+1}}{2}\right) - T^n x_n\right\| + \|T^n x_n - x_n\|$$

$$\leq a_n \|f(x_n) - T^n x_n\| + b_n \|x_n - T^n x_n\|$$
$$+ c_n \left\|T^n\left(\frac{x_n + x_{n+1}}{2}\right) - T^n x_n\right\| + \|T^n x_n - x_n\|$$

$$\leq a_n \|f(x_n) - T^n x_n\| + (1 + b_n)\|x_n - T^n x_n\| + c_n k_n \left\|\frac{x_n + x_{n+1}}{2} - x_n\right\|$$

$$= a_n \|f(x_n) - T^n x_n\| + (1 + b_n)\|x_n - T^n x_n\| + \frac{c_n k_n}{2}\|x_{n+1} - x_n\|. \tag{43}$$

Let $M > 0$ be a constant such that

$$M \geq \sup\{\|f(x_n) - T^n x_n\|, n \in \mathbb{N}\}. \tag{44}$$

Then,

$$\left(1 - \frac{c_n k_n}{2}\right)\|x_{n+1} - x_n\| \leq a_n M + (1 + b_n)\|x_n - T^n x_n\|,$$

$$\frac{2 - c_n k_n}{2}\|x_{n+1} - x_n\| \leq a_n M + (1 + b_n)\|x_n - T^n x_n\|,$$

$$\frac{1 + (1-\varepsilon)a_n + b_n}{2}\|x_{n+1} - x_n\| \leq a_n M + (1 + b_n)\|x_n - T^n x_n\|,$$

$$\|x_{n+1} - x_n\| \leq \frac{2a_n}{1 + (1-\varepsilon)a_n + b_n}M + \frac{2(1+b_n)}{1 + (1-\varepsilon)a_n + b_n}\|x_n - T^n x_n\|$$

$$2a_n M + 4\|x_n - T^n x_n\|. \tag{45}$$

Using conditions (i) and (iv),
$$\lim_{n \to \infty} \|x_{n+1} - x_n\| = 0. \tag{46}$$

Next, to evaluate $\lim_{n \to \infty} \|x_{n+1} - T^n x_n\|$,

$$\begin{aligned}
\|x_{n+1} - T^n x_n\| &= \left\|a_n f(x_n) + b_n x_n + c_n T^n\left(\frac{x_n + x_{n+1}}{2}\right) - T^n x_n\right\| \\
&\leq a_n \|f(x_n) - T^n x_n\| + b_n \|x_n - T^n x_n\| + c_n \left\|T^n\left(\frac{x_n + x_{n+1}}{2}\right) - T^n x_n\right\| \\
&\leq a_n \|f(x_n) - T^n x_n\| + b_n \|x_n - T^n x_n\| + c_n k_n \left\|\frac{x_n + x_{n+1}}{2} - x_n\right\| \\
&= a_n \|f(x_n) - T^n x_n\| + b_n \|x_n - T^n x_n\| + \frac{c_n k_n}{2} \|x_{n+1} - x_n\| \\
&\leq a_n M + b_n \|x_n - T^n x_n\| + \frac{c_n k_n}{2} \|x_{n+1} - x_n\| \longrightarrow 0,
\end{aligned}$$

as $n \longrightarrow \infty$. \hfill (47)

In a similar manner, one can apply (46) and (47) to obtain that

$$\begin{aligned}
\|x_n - T^{n-1} x_n\| &= \|x_n - T^{n-1} x_{n-1} + T^{n-1} x_{n-1} - T^{n-1} x_n\| \\
&\leq \|x_n - T^{n-1} x_{n-1}\| + \|T^{n-1} x_{n-1} - T^{n-1} x_n\| \\
&\leq \|x_n - T^{n-1} x_{n-1}\| + k_{n-1} \|x_n - x_{n-1}\| \longrightarrow 0,
\end{aligned}$$

as $n \longrightarrow \infty$. \hfill (48)

Consequently, from (46)–(48), it can be deduced that
$$\lim_{n \to \infty} \|x_n - T x_n\| = 0. \tag{49}$$

*Step 2.* It is claimed that
$$\limsup_{n \to \infty} \langle f(p) - p, J(x_n - p) \rangle \leq 0, \tag{50}$$

where $p$ is a fixed point of $T$.

Simply by applying Lemma 7, it is obtained that
$$\limsup_{n \to \infty} \langle f(p) - p, J(x_n - p) \rangle \leq 0. \tag{51}$$

Alternatively, define a function $\phi: K \longrightarrow \mathbb{R}^+$ by
$$\phi(z) = \mu_n \|x_n - z\|^2, \tag{52}$$

for all $z \in K$. Observe that $\phi$ is continuous, convex, and $\phi(z) \longrightarrow \infty$ as $\|z\| \longrightarrow \infty$. Since $E$ is reflexive, $\phi$ attains its infimum over $K$. Let $z_0 \in K$ such that $\phi(z_0) = \min_{z \in K} \phi(z)$ and let
$$G = \left\{x \in K: \phi(x) = \min_{z \in K} \phi(z)\right\}. \tag{53}$$

Then, $G$ is nonempty because $z_0 \in G$ and it is also closed, convex, and bounded (for more details, see e.g., [15, 19, 25]). Moreover, $T$ has a fixed point in $G$ by Lemma 6. Indeed, let $p \in G$ and $\omega_w(p)$ be the weak $w$-limit set of $T$ at $p$, that is, the set $\{x \in K: x: = w - \lim_j T^{n_j} p \text{ for some } n_j \longrightarrow \infty\}$. Then, from the weak-lower-semicontinuity of $\phi$ and the fact that $\lim_{n \to \infty} \|x_n - T^n x_n\| = 0$ (condition (iv)), the following estimate is obtained:

$$\begin{aligned}
\phi(x) &\leq \liminf_{j \to \infty} \phi(T^{n_j} p) \leq \limsup_{n \to \infty} \phi(T^n p) \\
&= \limsup_{n \to \infty} \mu_n \|x_n - T^n p\|^2 \\
&\leq \limsup_{n \to \infty} \mu_n \left(\|x_n - T^n x_n\| + \|T^n x_n - T^n p\|\right)^2 \\
&\leq \limsup_{n \to \infty} \mu_n \|T^n x_n - T^n p\|^2 \\
&\leq \limsup_{n \to \infty} \mu_n k_n^2 \|x_n - p\|^2 \\
&\leq \mu_n \|x_n - p\|^2 = \phi(p).
\end{aligned} \tag{54}$$

This shows that $G$ satisfies the property (P). By applying Lemma 6, $x = Tx$, it follows that $G \cap F(T) \neq \emptyset$. Suppose that $p \in G \cap F(T)$. Then, Lemma 1 gives that
$$\mu_n \langle x - p, J(x_n - p) \rangle \leq 0, \tag{55}$$

for all $x \in K$. In particular,
$$\mu_n \langle f(p) - p, J(x_n - p) \rangle \leq 0. \tag{56}$$

Since $\|x_{n+1} - x_n\| \longrightarrow 0$ as $n \longrightarrow \infty$ (46), and the duality mapping $J$ is norm to weak-star uniformly continuous on bounded set (earlier mentioned in Section 2); this gives
$$\lim_{n \to \infty} \left(\langle f(p) - p, J(x_{n+1} - p) \rangle - \langle f(p) - p, J(x_n - p) \rangle\right) = 0. \tag{57}$$

Therefore, the sequence $\{\langle f(p) - p, J(x_n - p) \rangle\}$ satisfies the conditions of Lemma 5. Thus,
$$\limsup_{n \to \infty} \langle f(p) - p, J(x_n - p) \rangle \leq 0. \tag{58}$$

*Step 3.* Lastly, it is shown that $x_n \longrightarrow p$:

$$\begin{aligned}\|x_{n+1} - p\|^2 &= a_n \langle f(x_n) - p, J(x_{n+1} - p) \rangle + b_n \langle x_n - p, J(x_{n+1} - p) \rangle \\
&\quad + c_n \langle T^n \left( \frac{x_n + x_{n+1}}{2} \right) - p, J(x_{n+1} - p) \rangle \\
&= a_n \langle f(x_n) - f(p), J(x_{n+1} - p) \rangle + b_n \langle x_n - p, J(x_{n+1} - p) \rangle \\
&\quad + c_n \langle T^n \left( \frac{x_n + x_{n+1}}{2} \right) - p, J(x_{n+1} - p) \rangle \\
&\leq \alpha a_n \|x_n - p\| \|x_{n+1} - p\| + b_n \|x_n - p\| \|x_{n+1} - p\| \\
&\quad + c_n k_n \left\| \frac{x_n + x_{n+1}}{2} - p \right\| \|x_{n+1} - p\| + a_n \langle f(p) - p, J(x_{n+1} - p) \rangle \\
&\leq \alpha a_n \|x_n - p\| \|x_{n+1} - p\| + b_n \|x_n - p\| \|x_{n+1} - p\| \\
&\quad + \frac{c_n k_n}{2} (\|x_n - p\| + \|x_{n+1} - p\|) \|x_{n+1} - p\| + a_n \langle f(p) - p, J(x_{n+1} - p) \rangle \\
&= \left( \alpha a_n + b_n + \frac{c_n k_n}{2} \right) \|x_n - p\| \|x_{n+1} - p\| + \frac{c_n k_n}{2} \|x_{n+1} - p\|^2 \\
&\quad + a_n \langle f(p) - p, J(x_{n+1} - p) \rangle \\
&= \frac{2\alpha a_n + 2b_n + c_n k_n}{2} \|x_n - p\| \|x_{n+1} - p\| + \frac{c_n k_n}{2} \|x_{n+1} - p\|^2 \\
&\quad + a_n \langle f(p) - p, J(x_{n+1} - p) \rangle \\
&\leq \frac{2\alpha a_n + 2b_n + c_n k_n}{4} \left( \|x_n - p\|^2 + \|x_{n+1} - p\|^2 \right) + \frac{c_n k_n}{2} \|x_{n+1} - p\|^2 \\
&\quad + a_n \langle f(p) - p, J(x_{n+1} - p) \rangle \\
&= \frac{2\alpha a_n + 2b_n + 3c_n k_n}{4} \|x_{n+1} - p\|^2 + \frac{2\alpha a_n + 2b_n + c_n k_n}{4} \|x_n - p\|^2 \\
&\quad + a_n \langle f(p) - p, J(x_{n+1} - p) \rangle.
\end{aligned} \tag{59}$$

Therefore,

$$\left(1 - \frac{2\alpha a_n + 2b_n + 3c_n k_n}{4}\right) \|x_{n+1} - p\|^2 \leq \frac{2\alpha a_n + 2b_n + c_n k_n}{4} \|x_n - p\|^2$$
$$+ a_n \langle f(p) - p, J(x_{n+1} - p) \rangle,$$
$$\frac{4 - (2\alpha a_n + 2b_n + 3c_n k_n)}{4} \|x_{n+1} - p\|^2 \leq \frac{2\alpha a_n + 2b_n + c_n k_n}{4} \|x_n - p\|^2$$
$$+ a_n \langle f(p) - p, J(x_{n+1} - p) \rangle,$$
$$\|x_{n+1} - p\|^2 \leq \frac{2\alpha a_n + 2b_n + c_n k_n}{4 - (2\alpha a_n + 2b_n + 3c_n k_n)} \|x_n - p\|^2$$
$$+ \frac{4a_n}{4 - (2\alpha a_n + 2b_n + 3c_n k_n)} \langle f(p) - p, J(x_{n+1} - p) \rangle. \tag{60}$$

Observe that

$$\begin{aligned}
4 - (2\alpha a_n + 2b_n + 3c_n k_n) &= 4 - 2\alpha a_n - 2b_n - 3k_n (1 - a_n - b_n) \\
&= 4 - 3k_n - 2\alpha a_n + 3k_n a_n - 2b_n + 3k_n b_n \\
&= 1 - 3(k_n - 1) - 2\alpha a_n + 3k_n a_n - 2b_n + 3k_n b_n \\
&\geq 1 - 3\varepsilon a_n - 2\alpha a_n + 3k_n a_n - 2b_n + 3k_n b_n \\
&= 1 + (3k_n - 2\alpha - 3\varepsilon)a_n + (3k_n - 2)b_n \\
&\geq 1 + (3 - 2\alpha - 3\varepsilon)a_n + (3 - 2)b_n \\
&= 1 + (3 - 2\alpha - 3\varepsilon)a_n + b_n,
\end{aligned}$$

for sufficiently large $n$. $\tag{61}$

Consequently,

$$\|x_{n+1} - p\|^2 \leq \frac{1 - (1 - 2\alpha - \varepsilon)a_n + b_n}{1 + (3 - 2\alpha - 3\varepsilon)a_n + b_n}\|x_n - p\|^2$$

$$+ \frac{4a_n}{1 + (3 - 2\alpha - 3\varepsilon)a_n + b_n}\langle f(p) - p, J(x_{n+1} - p)\rangle$$

$$= \left(1 - \frac{4(1 - \alpha - \varepsilon)a_n}{1 + (3 - 2\alpha - 3\varepsilon)a_n + b_n}\right)\|x_n - p\|^2$$

$$+ \frac{4(1 - \alpha - \varepsilon)a_n}{1 + (3 - 2\alpha - 3\varepsilon)a_n + b_n} \frac{1}{(1 - \alpha - \varepsilon)}\langle f(p) - p, J(x_{n+1} - p)\rangle. \quad (62)$$

By applying Lemma 2 to (58) and (62), it is deduced that $x_n \longrightarrow p$ as $n \longrightarrow \infty$,

where, respectively, $\gamma_n := 4(1 - \alpha - \varepsilon)a_n/1 + (3 - 2\alpha - 3\varepsilon)a_n + b_n$ and $\beta_n := 1/1 - \alpha - \varepsilon\langle f(p) - p, J(x_{n+1} - p)\rangle$.

*Remark 1.* It is generally known that a Hilbert space is a subclass of uniformly convex Banach spaces. Moreover, it is obvious that every nonexpansive mapping is an asymptotically nonexpansive mapping. Consequently, Theorem 1 is an improvement and generalization of the main results of Zhao et al. [26], Aibinu et al. [10], Luo et al. [27], and Yao et al. [7] as stated below.

**Corollary 1** (see [26]). *Let K be a nonempty closed and convex subset of a Hilbert space H, $T: K \longrightarrow K$ be an asymptotically nonexpansive mapping with a sequence $\{k_n\} \subset [1, +\infty)$, $\lim_{n \to \infty} k_n = 1$ and $F(T) \neq \emptyset$. Let f be a contraction on K with coefficient $c \in [0, 1)$. For an arbitrary initial point $x_1 \in K$, let $\{x_n\}$ be the sequence generated by*

$$x_{n+1} = a_n f(x_n) + (1 - a_n)T^n\left(\frac{x_n + x_{n+1}}{2}\right), \quad n \in \mathbb{N}. \quad (63)$$

*Then, the sequence $\{x_n\}$ converges strongly to $p = P_{F(T)}f(p)$ which solves the following variational inequality (8).*

*Proof.* Take $E$ to be a Hilbert space in Theorem 1. Moreover, since $a_n + b_n + c_n = 1, \forall n \in \mathbb{N}$, by taking $\{b_n\} = \{0\}$ in (10), we have $c_n = 1 - a_n$. Therefore, we obtain (63) from (10) by taking $b_n = 0$. Hence, the desired result follows from Theorem 1. □

**Corollary 2** (see [10]). *Let E be a uniformly smooth Banach space and K be a nonempty closed convex subset of E. Let $T: K \longrightarrow K$ be a nonexpansive mapping with $F(T) \neq \emptyset$ and $f: K \longrightarrow K$ be a c-contraction. Suppose $\{a_n\}$ satisfies (i) and (ii) and $\{b_n\}$ satisfies*

($C_1$) $0 < \liminf_{n \to \infty} \beta_n \leq \limsup_{n \to \infty} \beta_n < 1$

($C_2$) $\lim_{n \to \infty}|\beta_{n+1} - \beta_n| = 0$

*For an arbitrary $x_1 \in K$, define the iterative sequence $\{x_n\}$ by (9). Then, as $n \longrightarrow \infty$, the sequence $\{x_n\}$ converges in norm to a fixed point p of T, where p is the unique solution in $F(T)$ to the variational inequality (8).*

*Proof.* Every nonexpansive mapping is an asymptotically nonexpansive mapping and all uniformly convex or uniformly smooth Banach spaces have uniform normal structure (see [14, 15]). Therefore, it suffices to show that $\lim_{n \to \infty}\|x_n - Tx_n\| = 0$:

$$\|x_n - Tx_n\| \leq \|x_n - x_{n+1}\| + \|x_{n+1} - Tx_n\|$$

$$= \|x_n - x_{n+1}\| + a_n\|f(x_n) - Tx_n\| + b_n\|x_n - Tx_n\|$$

$$+ c_n\left\|T\left(\frac{x_n + x_{n+1}}{2}\right) - Tx_n\right\|$$

$$= \|x_n - x_{n+1}\| + a_n\|f(x_n) - Tx_n\| + b_n\|x_n - Tx_n\| + c_n\left\|\frac{x_n + x_{n+1}}{2} - x_n\right\|$$

$$= \|x_n - x_{n+1}\| + a_n\|f(x_n) - Tx_n\| + b_n\|x_n - Tx_n\| + \frac{c_n}{2}\|x_n - x_{n+1}\|. \quad (64)$$

Let $M$ be a constant such that $M > \sup_{n \in N}\|f(x_n) - Tx_n\|$ and $b = \sup_n (1 - b_n) > 0$ since $b_n < 1$. Then,

$$\|x_n - Tx_n\| \leq \frac{2 + c_n}{2(1 - b_n)}\|x_{n+1} - x_n\| + \frac{a_n}{1 - b_n}M$$

$$\leq \frac{2 + c_n}{2b}\|x_{n+1} - x_n\| + \frac{a_n}{b}M \longrightarrow 0, \quad \text{as } n \longrightarrow \infty. \quad (65)$$

This completes the proof. □

**Corollary 3** (see [27]). *Let K be a closed convex subset of a uniformly smooth Banach space E. Let $T: K \longrightarrow K$ be a nonexpansive mapping with $F(T) \neq \emptyset$ and $f: K \longrightarrow K$ a contraction with coefficient $c \in [0, 1)$. Let $\{x_n\}$ be a sequence generated by the viscosity implicit midpoint rule (7), where $\{a_n\}$ is a sequence in $(0, 1)$ such that it satisfies conditions (i) and (ii) and either $(C_3) \sum_{n=1}^{\infty} |a_n - a_{n-1}| < +\infty$ or $\lim_{n \to \infty}(a_n/a_{n-1}) = 1$. Then, $\{x_n\}$ converges strongly to a fixed point p of T, which also solves the variational inequality (8).*

*Proof.* Take $b_n = 0$ in equation (9). The desired result follows from Corollary (28). □

**Corollary 4** (see [7]). *Let K be a nonempty closed convex subset of a real Hilbert space H. Let $T: K \longrightarrow K$ be a nonexpansive mapping with $F(T) \neq \emptyset$. Suppose $f: K \longrightarrow K$ is a c-contraction. For given $x_0 \in K$ arbitrarily, let the sequence $\{x_n\}$ be generated by (9). Assume that $\{a_n\}$ satisfies (i) and (ii) and $\{b_n\}$ satisfies*

($C_1$) $0 < \liminf_{n \to \infty} b_n \leq \limsup_{n \to \infty} b_n < 1$

($C_4$) $\lim_{n \to \infty}(b_{n+1} - b_n) = 0$

*Then, the sequence $\{x_n\}$ generated by (9) converges strongly to $p = P_{F(T)}Q(p)$.*

*Proof.* Observe that a Hilbert space is a subclass of a uniformly smooth Banach space. By taking $E$ to be a Hilbert space in Corollary 2, the result of Yao et al. in [7] is obtained. □

Table 1: Numerical values for the $n$th iterations as $\{x_n\}$ converges to a fixed point.

| Iteration ($n$) | $\{x_{n+1}\}$ | | |
|---|---|---|---|
| | ($p = (1,-1), x_1 = (0, 1/3)$) | ($p = (0,0), x_1 = (1/2, 1)$) | ($p = (-1, 1), x_1 = (-2, 1)$) |
| 1 | 0.5152 | 0.2236 | 0.6485 |
| 2 | 0.0008 | 0.0003 | 0.0010 |
| 3 | 0.0003 | 0.0001 | 0.0004 |
| 4 | 0.0002 | 0.0001 | 0.0002 |
| 5 | 0.0001 | 0.0001 | 0.0002 |
| 6 | 0.0001 | 0.0000 | 0.0001 |
| 7 | 0.0001 | 0.0000 | 0.0001 |
| 8 | 0.0001 | 0.0000 | 0.0001 |
| 9 | 0.0001 | 0.0000 | 0.0001 |
| 10 | 0.0000 | 0.0000 | 0.0001 |
| 11 | 0.0000 | 0.0000 | 0.0000 |
| 12 | 0.0000 | 0.0000 | 0.0000 |
| 13 | 0.0000 | 0.0000 | 0.0000 |
| 14 | 0.0000 | 0.0000 | 0.0000 |
| 15 | 0.0000 | 0.0000 | 0.0000 |
| 16 | 0.0000 | 0.0000 | 0.0000 |
| 17 | 0.0000 | 0.0000 | 0.0000 |
| 18 | 0.0000 | 0.0000 | 0.0000 |
| 19 | 0.0000 | 0.0000 | 0.0000 |
| 20 | 0.0000 | 0.0000 | 0.0000 |

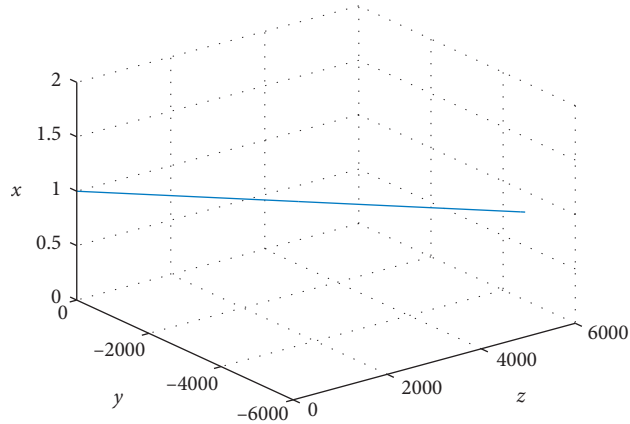

Figure 1: Convergence of $\{x_n\}$ to the fixed point $p = (1, -1)$ with $x_1 = (0, 1/3)$.

## 4. Numerical Example

Numerical example is given in this section to illustrate the convergence of the sequence of iteration in the main theorem.

*Example 1.* Let $K$ be an orthogonal subspace of $\mathbb{R}^2$ with the norm $\|u\| = \sqrt{\sum_{i=1}^{2} u_i^2}$ and the inner product $\langle u, v \rangle = \sum_{i=1}^{2} u_i v_i$ for $u = (u_1, u_2)$ and $v = (v_1, v_2)$. Define a mapping $T: K \longrightarrow K$ by

$$Tu = \begin{cases} (u_1, u_2), & \text{if } \Pi_{i=1}^{2} u_i < 0, \\ (-u_1, -u_2), & \text{if } \Pi_{i=1}^{2} u_i \geq 0, \end{cases} \quad (66)$$

for each $u = (u_1, u_2) \in K$. Then,

$$\|T^n u - T^n v\|^2 = \|u - (-1)^n v\|^2 = \|u\|^2 + \|v\|^2 = \|u - v\|^2, \quad (67)$$

or

$$\|T^n u - T^n v\|^2 = \|(-1)^n u - v\|^2 = \|u\|^2 + \|v\|^2 = \|u - v\|^2, \quad (68)$$

for any $u, v \in K$ (Zhang et al. [28]). Take $\{k_n\} \subset [1, \infty)$ to be $k_n := \{1 + (1/2^n)\}$ and observe that $\|T^n u - T^n v\| \leq k_n \|u - v\|$ for any $u, v \in K$. Thus, $T$ is an asymptotically nonexpansive mapping on $K$ and $F(T) = \{(0,0) \cup (u_1, u_2)\}$. We define $f: \mathbb{R} \longrightarrow \mathbb{R}$ by $f(x) = (1/2)x$. Let $a_n = 1/n, b_n = n - 1/n(n+1)$, and $c_n = n - 1/n + 1$. It is obvious that the sequences $\{a_n\}, \{b_n\}$, and $\{c_n\}$ satisfy conditions (i)–(iii) in Assumption 1, and according to (67) or (68), $\lim_{n \longrightarrow \infty} \|u_n - T^n u_n\| = \|u - u\| = 0$. The table and figures are given to display the numerical results and convergence of our sequence of iteration to fixed points of $T$. The computations are carried out by using MATLAB. Table 1 exhibits numerical values for the nth iterations as $\{x_n\}$ converges to a fixed point. Figure 1 displays the convergence of $\{x_n\}$ to the fixed point $p = (1, -1)$ with $x_1 = (0, 1/3)$. Convergence of $\{x_n\}$ to the fixed point $p = (1, -1)$ with $x_1 = (0, 1/3)$ is shown in Figure 2. Graph of $\|x_{n+1} - p\|$ as $\{x_n\}$ converges to the fixed point $p = (-1, 1)$ with $x_1 = (-2, 1)$ is presented in Figure 3.

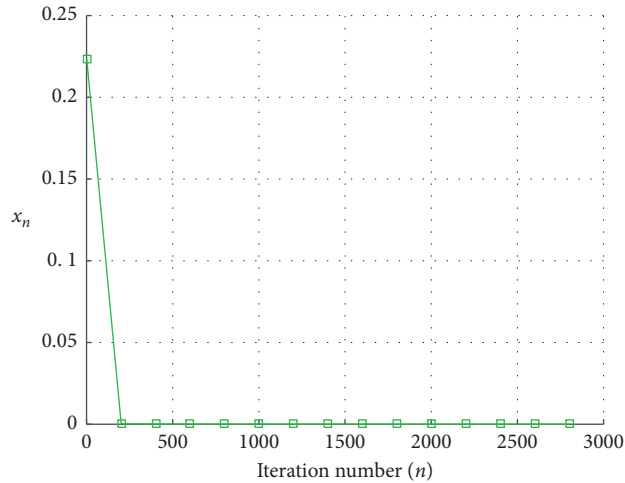

Figure 2: Convergence of $\{x_n\}$ to the fixed point $p = (0, 0)$ with $x_1 = (1/2, 1)$.

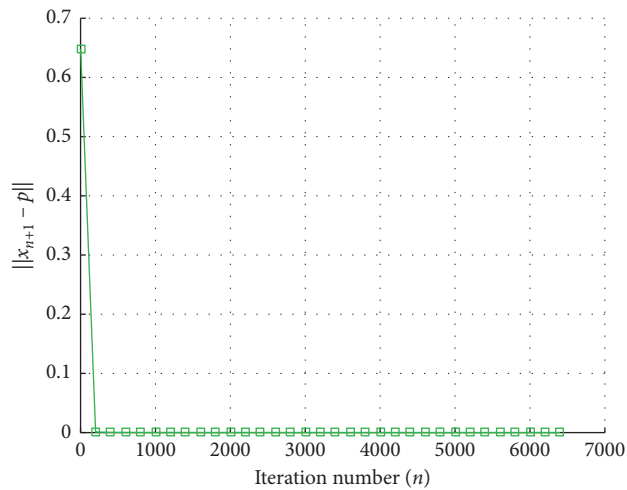

Figure 3: Graph of $\|x_{n+1} - p\|$ as $\{x_n\}$ converges to the fixed point $p = (-1, 1)$ with $x_1 = (-2, 1)$.

## 5. Conclusion

This paper considered the implicit midpoint rule of asymptotically nonexpansive mappings using the viscosity technique. The necessary conditions for the convergence of the class of asymptotically nonexpansive mappings are established by using this well-known iterative algorithm which plays important roles in the computation of fixed points of nonlinear mappings. Moreover, the previous results have been extended from the Hilbert space to a Banach space with a uniformly Gateaux differentiable norm possessing uniform normal structure. A numerical example is given to show the convergence of the sequence of iteration in the main theorem. The result presented in this paper plays crucial theoretical base of the numerical analysis for nonlinear problems which are often encountered in physical and biological sciences. The future direction for this paper is to determine for the class of asymptotically nonexpansive mappings and the conditions for the convergence of Jungck Noor iteration with $s$-convexity, defined by Kang et al. [2] and to compare the rate of convergence with the implicit midpoint procedure.

## Data Availability

The data used to support the findings of this study are available from the corresponding author upon request.

## Conflicts of Interest

The authors declare no conflicts of interest.

## Authors' Contributions

All authors contributed significantly in writing this article. All authors read and approved the final manuscript.